\newcommand{\force}{\Vdash} 
\newtheorem{theorem}{Theorem}
\newtheorem{lemma}[theorem]{Lemma}

\newtheorem{proposition}[theorem]{Proposition}
\newtheorem{definition}[theorem]{Definition}

\newtheorem{corollary}[theorem]{Corollary}
\newtheorem{remark}[theorem]{Remark}
\documentclass[12pt]{article}
\def\ma{{\mathfrak A}}
\def\mb{{\mathfrak B}}
\def\proof{\noindent {\bf Proof.}\hspace{2mm}}
\def\qed{$\Box$\medbreak }
\def\n{\nu}\def\g{\gamma}\def\d{\delta}\def\k{\kappa}
\def\raj{\restriction}

\def\cf{{\rm cof}}
\def\coinit{{\rm coinit}}
\def\iso{\cong}
\def\rank{{\rm rnk}}%
\def\emp{{\forall}}
\def\non{{\exists}}
\usepackage{latexsym,amssymb}
\title{More on the Ehrenfencht-Fra\"\i ss\'e game of length $\omega_1$.}
\author{
Tapani Hyttinen
\thanks{Partially supported by the Academy of
Finland grant \#40734.}\\ Department of  Mathematics\\ University
of Helsinki\\ Helsinki, Finland \and Saharon Shelah
\thanks{Research partially supported by the United States-Israel
        Binational Science Foundation. Publication number [HShV:776] }\\
Institute of Mathematics\\ Hebrew University\\ Jerusalem, Israel\\
\and Jouko V\"a\"an\"anen
\thanks{Partially supported by the Academy of
Finland grant \#40734.}\\ Department of  Mathematics\\ University
of Helsinki\\ Helsinki, Finland}

\begin{document}
\maketitle
\def\qed{$\Box$\medbreak }
\def\I{\forall}
\def\II{\exists}
\def\la{\langle}
\def\ra{\rangle}
\newcommand{\EFG}[1]{{\rm EFG}_{#1}}
\def\BM{{\rm PDG}}
\def\A{{\cal A}}
\def\B{{\cal B}}
\def\C{{\cal C}}
\def\D{{\cal D}}
\def\E{{\cal E}}
\def\F{{\cal F}}\def\H{{\cal H}}
\def\P{{\cal P}}
\def\PP{{\mathbb P}}
\def\QQ{{\mathbb Q}}
\def\L{{\cal L}}
\def\M{{\bf M}}
\def\K{{\bf K}}
\def\Q{{\mathbb Q}}
\def\Z{{\mathbb Z}}
\def\R{{\mathbb R}}
\def\J{{\cal J}}
 \def\I{{\cal I}}
 \def\G{{\cal G}}
 \def\ef{Ehrenfencht-Fra\"\i ss\'e }
{\def\sqr#1#2{{\vcenter{\hrule height.#2pt
      \hbox{\vrule width.#2pt height#1pt \kern#1pt
         \vrule width.#2pt}
       \hrule height.#2pt}}}
\def\e{\epsilon}
\def\k{\kappa}
\def\o{\omega}
\def\s{\sigma}
\def\d{\delta}
\def\v{\varphi}
\def\l{\lambda}
\def\r{\rho}
\def\a{\alpha}
\def\b{\beta}
\def\p{\pi}
\def\z{\zeta}
\def\t{\tau}
\def\g{\gamma}
\def\h{\theta}
\def\x{\xi}
\def\n{\eta}
\def\ga{\alpha}
\def\gb{\beta}
\def\gp{\pi}
\def\gz{\zeta}
\def\gt{\tau}
\def\gg{\gamma}
\def\go{\omega}
\def\goy{\omega_1}
\def\gk{\kappa}
\def\gl{\lambda}
\def\gs{\sigma}
\def\se{\subseteq}


\def\abs#1{\left\vert #1\right\vert}
\def\norm#1{\left\Vert #1\right\Vert}
\def\rem #1 #2. #3\par{\medbreak{\bf #1 #2.
\enspace}{#3}\par\medbreak}
\def\proof{{\bf Proof}.\enspace}
\def\sqr#1#2{{\vcenter{\hrule height.#2pt
      \hbox{\vrule width.#2pt height#1pt \kern#1pt
         \vrule width.#2pt}
       \hrule height.#2pt}}}
\def\eop{\mbox{\qed}}

\def\reals{\hbox{\rm I\hskip-0.14em R}}
\def\Nat{\hbox{\rm I\hskip-0.14em N}}
\def\N{\hbox{\rm I\hskip-0.14em N}}
\def\Re{\hbox{\rm Re}\;}
\def\Im{\hbox{\rm Im}\;}
\def\sgn{\hbox{\rm sgn}\;}
\def\smallRe{\hbox{\sevenrm Re}}
\def\smallIm{\hbox{\sevenrm Im}}
\def\sm{\ll}
\def\ol{\overline}
\def\tr{{\rm tr}}
\def\empty{\emptyset}
\def\at{{\rm at}}
\def\raj{\restriction}
\def\tar{\models}
\def\inv{^{-1}}
\def\sm{\ll}
\def\ol{\overline}
\def\mod{{\rm \ mod \ }}
\def\fork{\hbox{$\ni\hskip -2mm-$}}
\def\nfork{\hbox{$\ni\hskip -2mm\not-$}}
\def\da{\downarrow}
\def\Da{\Downarrow}
\def\kom{\mapsto}
\def\nda{\mathrel{\lower0pt\hbox to 3pt{\kern3pt$\not$\hss}\downarrow}}
\def\nDa{\mathrel{\lower0pt\hbox to 3pt{\kern3pt$\not$\hss}\Downarrow}}
\def\nbot{\mathrel{\lower0pt\hbox to 4pt{\kern3pt$\not$\hss}\bot}}
\def\ekom{\mathrel{\lower3pt\hbox to 0pt{\kern3pt$\sim$\hss}\mapsto}}
\def\do{\triangleright}
\def\edo{\trianglerighteq}
\def\ndo{\not\triangleright}
\def\T{\Theta}
\def\anR{\mathrel{\lower1pt\hbox to 2pt{\kern3pt$R$\hss}\not}}
\def\anoR{\mathrel{\lower1pt\hbox to 2pt{\kern3pt$\overline{R}$\hss}\not}}
\def\nR{\anR\ \ }
\def\noR{\anoR\ \ }
\def\anRm{\mathrel{\lower1pt\hbox to 2pt{\kern3pt$R^{-}$\hss}\not}}
\def\nRm{\anRm\ \ \ \ }
\def\ndda{\mathrel{\lower0pt\hbox to 1pt{\kern3pt$\not$\hss}\downdownarrows}}
\def\dda{\downdownarrows}
\def\warrow{\mathrel{\lower0pt\hbox to 1pt{\kern3pt$^{w}$\hss}\rightarrow}}
\def\bbox{\blacksquare}
\def\bdia{\blacklozenge}


This paper is a continuation of \cite{msv416}. Let \(\ma\) and
\(\mb\) be two first order structures of the same vocabulary
\(L\). We denote the domains of \(\ma\) and \(\mb\) by \(A\) and
\(B\) respectively. All vocabularies are assumed to be relational.
The {\em Ehrenfeucht-Fra\"\i ss\'e-game of length \(\gg\) of
\(\ma\) and \(\mb\)} denoted  by \(\EFG{\gg}(\ma,\mb)\) is defined
as follows: There are
 two players
called \(\forall\) and \(\exists\). First \(\forall\) plays
\(x_0\) and then \(\exists\) plays \(y_0\). After this
 \(\forall\) plays \(x_1\), and \(\exists\) plays \(y_1\), and
so on. If  \(\la(x_{\gb},y_{\gb}):\gb<\ga\ra\) has been played and
\(\ga<\gg\), then  \(\forall\) plays \(x_{\ga}\) after which
\(\exists\) plays \(y_{\ga}\). Eventually a sequence
\(\la(x_{\gb},y_{\gb}):\gb<\gg\ra\) has been played. The rules of
the game say that both players have to play elements of \(A\cup
B\). Moreover, if \(\forall\) plays his \(x_{\gb}\) in \(A\)
(\(B\)), then \(\exists\) has to play his \(y_{\gb}\) in \(B\)
(\(A\)). Thus the sequence \(\la(x_{\gb},y_{\gb}):\gb<\gg\ra\)
determines a relation \(\pi\se A\times B\). Player \(\exists\)
wins this round of the game if \(\pi\) is a partial isomorphism.
Otherwise \(\forall\) wins. The notion of winning strategy is
defined in the usual manner. The game
\(\EFG{\gamma}^\delta(\ma,\mb)\) is defined like
\(\EFG{\gamma}(\ma,\mb)\) except that the players play sequences
of length \(<\delta\) at a time. Thus \(\EFG{\gamma}(\ma,\mb)\) is
the same game as \(\EFG{\gamma}^2(\ma,\mb)\).

It was proved in \cite{msv416} that, assuming \(\square_{\goy}\),
there are models \(\ma\) and \(\mb\) of cardinality \(\aleph_2\)
such that the game \({\cal G}_{\goy}(\ma,\mb)\) is non-determined.
In this paper we weaken the assumption \(\square_{\goy}\), to
``\(\omega_2\) is not weakly compact in \(L\)''
(Corollary~\ref{wc}), but we can do this only if we assume CH. We
do not know if this is possible without CH. In the other
direction, it was proved in \cite{msv416} that if the
\(\goy\)-nonstationary ideal on \(\omega_2\) has a \(\gs\)-closed
dense subset, then the
 game \(\EFG{\goy}(\ma,\mb)\) is determined for
all \(\ma\) and \(\mb\) of cardinality \(\le\aleph_2\). The
assumption is equivconsistent with the existence of a measurable
cardinal. In this paper we weaken the assumption to a condition
which is consistent relative to the existence of a weakly compact
cardinal (Corollary~\ref{i}). Thus we establish:

\begin{theorem}\label{main1}
The following statements are equiconsistent relative to ZFC:
\begin{enumerate}
\item
There is a weakly compact cardinal.
\item
CH and \(EF_{\omega_1}(\ma,\mb)\) is determined for all models
\(\ma\) and \(\mb\) of cardinality \(\aleph_2\).
\end{enumerate}
\end{theorem}

In \cite{msv416} we proved in ZFC that there are structures
\(\ma\) and \(\mb\) of cardinality \(\aleph_3\) with one binary
predicate such that the game \(\EFG{\goy}(\ma,\mb)\) is
non-determined. We now improve this result under some cardinal
arithmetic assumptions. We prove:

\begin{theorem}\label{1}
Assume that $2^{\o}<2^{\o_{3}}$ and $T$ is a countable complete
first order theory. Suppose that one of (i)-(iii) below holds.
Then there are $\A ,\B\models T$ of power $\o_{3}$ such that for
all cardinals $1<\theta\le\o_{3}$, $EF^{\theta}_{\o_{1}}(\A ,\B )$
is non-determined.
\begin{description}
\item[(i)] $T$ is unstable.

\item[(ii)] $T$ is superstable with DOP or OTOP.

\item[(iii)] $T$ is stable and unsuperstable and $2^{\o}\le\o_{3}$.
\end{description}
\end{theorem}

This result complements the result in \cite{msv416}  that if \(T\)
is an $\go$-stable
 first order theory with NDOP, then
\(\EFG{\goy}(\ma,\mb)\) is determined for all models \(\ma\) of
\(T\) and all models \(\mb\). This is actually true under the
weaker assumption that \(T\) is superstable with NDOP and NOTOP.

{\bf Notation:} We follow Jech \cite{Jech} in set theoretic
notation. We use \(S^m_n\) to denote the set \(\{\ga<\go_m : \cf
(\ga)=\go_n\}\). Closed and unbounded sets are called cub sets. A
set  of ordinals is {\em \(\gl\)-closed} if it is closed under
supremums of ascending \(\gl\)-sequences \(\la \ga_i : i<\gl\ra\)
of its elements. A subset of a cardinal is {\em
\(\gl\)-stationary} if it meets every \(\gl\)-closed unbounded
subset of the cardinal.

\section{Getting a weakly compact cardinal}

In this section we show that if CH holds and
\(\EFG{\goy}(\ma,\mb)\) is determined for all models \(\ma\) and
\(\mb\) of cardinality \(\aleph_2\), then \(\go_2\) is weakly
compact in \(L\) (Corollary~\ref{wc}). We use the results from
\cite{mag} that if \(\go_2\) is not weakly compact in \(L\), then
there is a bistationary \(S\se S^2_0\) such that for all
\(\alpha<\omega_2\) either \(\alpha\cap S\) or \(\alpha\setminus
S\) is non-stationary.

If \(I\) is a linear order, we use \((I)^*\) to denote the reverse
order of \(I\). We call a sequence \(s=(s_\xi)_{\xi<\zeta}\) {\em
coinitial sequence of length \(\zeta\)} in \(I\), if it is
decreasing in \(I\) and has no lower bound in \(I\). The {\em
coinitiality} \(\coinit(I)\) of a linear order \(I\) is the
smallest length of a coinitial sequence in \(I\).

Let \(\theta=\omega+ ((\omega_1)^*+\omega)\cdot\omega_1\).

\begin{lemma}\label{dense}
There is a dense linear order \(I\) such that
\begin{description}
\item[(i)] \(|I|=\aleph_1\).
\item[(ii)] \(\coinit(I)=\aleph_0\)
\item[(iii)] \(I\cdot(\alpha+1)\cong I\) for all \(\alpha\le\omega_1\).
\item[(iv)] \(I\cong I\cdot\omega+I\cdot(\omega_1)^*\).
\item[(v)]\(I\cdot\theta+I\cong I\).
\end{description}
\end{lemma}

\proof This is like Lemma 4.7.16 in \cite{hhr}. If \(J_1\) and
\(J_2\) are linear orders, let \(H(J_1,J_2)\) be the set of
\(f:n_f\rightarrow J_1\cup J_2\), where \(n_f<\omega\) is even,
\(f(2i)\in J_1\) and \(f(2i+1)\in J_2\) for all \(i<n_f\). We can
make \(H(J_1,J_2)\) a linear order by ordering the functions
lexicographically, i.e.
\[f\le g\iff \exists m\le n_f
(\forall i<m(f(i)=g(i))\& (m<n_f\rightarrow f(m)<g(m))).\] Let
$I_0=H(\Q,\omega+(\omega_1)^*)$ and \(I_1=H(I_0,\omega_1)\). Thus
\(I_0\cong(1+I_0)\cdot (\omega+(\omega_1)^*)\cdot\Q\) and
\(I_1\cong(1+I_1)\cdot\omega_1\cdot I_0\). By using
\(\Q\cong\Q+1+\Q\), \(\omega=1+\omega\) and
\(\omega_1=1+\omega_1\), one gets easily the following, first for
\(I_0\), and then for \(I_1\):
\begin{equation}\label{7}\label{2}
I_0\cong I_0+1+I_0\hspace{2mm} ,\hspace{2mm} I_1\cong I_1+1+I_1.
\end{equation}

Let \(I\) be the set of \(f:\omega\rightarrow I_1\cup \theta\),
where \(f(2i)\in I_1\) and \(f(2i+1)\in \theta\) for all
\(i<\omega\) ordered lexicographically. Thus \(I\cong
I\cdot\theta\cdot I_1\). In fact, \(I\) is of the form
\(J\cdot\Q\), so (ii) is true. By (\ref{2}) and \(\theta\cong
1+\theta\) one gets immediately (v). As \(I\cong I\cdot\theta\cdot
(1+I_1)\cdot\omega_1\cdot I_0\), we get from (v) easily (iii) for
\(\alpha=\omega_1\). From this and \(\alpha+\omega_1= \omega_1\)
we get immediately (iii) for \(\alpha<\omega_1\). Note that
\(\theta\cong\omega+(\omega_1)^*+\theta\). If we combine this with
\(I\cong I\cdot\theta\cdot I_1\) and
\((\omega_1)^*\cong(\omega_1)^*+1 \), we get (iv).

As to (i), we only have \(|I|=2^\omega\). We use this lemma in a
context where CH is assumed, so we could simply assume it here.
But actually the lemma is true without CH, as we can construct
\(I\) in \(L\). Then \(|I|=\aleph_1\). Note that our \(I_0\) and
\(I_1\) are in \(L\), and the only property of \(\omega_1\) that
we used was that it is a limit ordinal. \qed

\begin{definition}
Suppose \(S\subseteq S^2_0\). We define
\[\Phi(S)=\sum_{i<\omega_2}\eta_i,\]
where
\[\eta_i=\left\{
\begin{array}{ll}
I\cdot(\omega_1)^*,&\mbox{ if }i\in S\\ I,&\mbox{ if }i\notin S.
\end{array}\right.\]
Let \(\Phi_{\alpha,\beta}(S)\) be the suborder \(\sum_{\alpha\le
i<\beta}\eta_i\) of \(\Phi(S)\). The {\em rank} of \(x\in\Phi(S)\)
is the least \(\alpha\) such that \(x\in
\Phi_{\alpha,\alpha+1}(S)\). We denote this \(\alpha\) by
\(\rank(\Phi(S),x)\).\end{definition}

\begin{lemma}\label{isom1}
Assume \(S\subseteq S^2_0\) is such that there is no \(\alpha\in
S^2_1\) with both \(S\cap\alpha\) and \((S\cap S^2_0)\setminus S\)
stationary. Then
\[\Phi_{\alpha,\beta+1}(S)\cong I\]
whenever \(\alpha<\beta<\omega_2\) and \(\alpha\notin S\).
\end{lemma}
{\def\F{\Phi_{\alpha,\beta+1}(S)} \proof This is like Lemma 4.7.19
in \cite{hhr}. We use Lemma~\ref{dense} and induction on
\(\beta\).

Let us first assume \(\beta\notin S\). If \(\beta\) is a successor
ordinal, then \(\F\iso I+I=I\) by (iii). If \(\beta\) has
cofinality \(\omega\), then \(\F\iso I\cdot\omega+I\iso I\). If
\(\beta\) has cofinality \(\omega_1\) and \(\beta\cap S\) is
non-stationary, then \(I\iso I\cdot\omega_1+I\iso I\). Finally, if
\(\beta\) has cofinality \(\omega_1\) and \(\beta\setminus S\) is
non-stationary, then \(I\iso I\cdot\theta+I\iso I\), by (v).

Let us then assume \(\beta\in S\). Thus \(\beta\) has cofinality
\(\omega\). Therefore \(\F\iso I\cdot\omega+
I\cdot(\omega_1)^*\iso I\), by (iv). \qed}

\begin{lemma}\label{isom2}
Assume \(S\subseteq S^2_0\) is such that there is no \(\alpha\in
S^2_1\) with both \(S\cap\alpha\) and \((S\cap S^2_0)\setminus S\)
stationary. Then \(\Phi_{0,\alpha}(S)\cong
\Phi_{0,\alpha}(\emptyset)\) whenever \(\alpha\in S^2_1\) and
\(S\cap\alpha\) is not stationary.
\end{lemma}

\proof Let \((\alpha_\xi)_{\xi<\omega_1}\) by a continuously
increasing cofinal sequence in \(\alpha\) such that
\(\alpha_\xi\notin S\) for all \(\xi<\omega_1\). By
Lemma~\ref{isom1} there is an isomorphism
\[f_\xi:\Phi_{\alpha_\xi,\alpha_{\xi+1}+1}(S)\rightarrow
\Phi_{\alpha_\xi,\alpha_{\xi+1}+1}(\emptyset).\] Let
\(f=\cup_{\xi<\omega_1} f_\xi\). This is the required isomorphism.
\qed

\begin{proposition}\label{nondetlo}
Assume CH and that there is \(S\subseteq S^2_0\)  such that both
\(S\) and \(S^2_0\setminus S\) are stationary but there is no
\(\alpha\in S^2_1\) with both \(S\cap\alpha\) and \((S\cap
S^2_0)\setminus S\) stationary. Then there are models \(\ma\) and
\(\mb\) of cardinality \(\aleph_2\) such that
\(EF_{\omega_1}(\ma,\mb)\) is non-determined.
\end{proposition}

\proof We may assume, that \(\{\alpha\in S^2_1:\alpha\cap S \mbox{
is non-stationary }\}\) is stationary, for otherwise we work with
\(S'=S^2_0\setminus S\). Let \(\ma=\Phi(S)\) and
\(\mb=\Phi(\emptyset)\). We first show that \(\exists\) cannot
have a winning strategy in \(EF_{\omega+\omega+1}(\ma,\mb)\).
Suppose \(\tau\) is a strategy of \(\exists\). Let \(C\) be the
cub of ordinals \(\alpha<\omega_2\) such that if during the first
\(\omega\) rounds of the game, \(\forall\) plays elements of the
models of rank \(<\alpha\), then so does \(\exists\) following
\(\tau\). Let \(\delta\in C\cap S\). Let \((\delta_n)_{n<\omega}\)
be an increasing cofinal sequence in \(\delta\). Now we let
\(\forall\) play against \(\tau\) as follows: On round number
\(n<\omega\) we let \(\forall\) play some element of \(\ma\), if
\(n\) is even, and of \(\mb\), if \(n\) is odd, of rank
\(\delta_n\). During rounds \(\omega+n\),
\(n<\omega\), we let \(\forall\) play a coinitial sequence
of length \(\omega\) in
\(\Phi_{\delta,\delta+1}(\emptyset)\subseteq\ma\). As
\(\coinit(\Phi_{\delta,\delta+1}(S))=\omega_1\), 
the game is lost for \(\exists\). So
\(\tau\) could not be a winning strategy.

Suppose then \(\rho\) is a strategy of \(\forall\). We show that
this cannot be a winning strategy. By CH we have an
\(\omega_1\)-cub set \(D\)
of ordinals \(\delta<\omega_2\) such that if \(\exists\) plays
only elements of rank \(<\delta\), then \(\rho\) directs
\(\forall\) to play also elements of rank \(<\delta\) only. Let
\(\delta\in D\cap S^2_1\) such that \(\delta\cap S\) is
non-stationary. By Lemma~\ref{isom2} there is an isomorphism
\(f:\Phi_{0,\alpha}(S)\rightarrow \Phi_{0,\alpha}(\emptyset)\).
Now \(\exists\) can beat \(\rho\) by using \(f\). \qed

\begin{corollary}\label{wc}
If CH holds and \(\EFG{\goy}(\ma,\mb)\) is determined for all
models \(\ma\) and \(\mb\) of cardinality \(\aleph_2\), then
\(\go_2\) is weakly compact in \(L\).
\end{corollary}


\section{Getting determinacy from a weakly compact cardinal}


In this section we show that if \(\kappa\) is weakly compact, then
there is a forcing extension in which the game
\(\EFG{\goy}(\ma,\mb)\) is determined for all \(\ma\) and \(\mb\)
of cardinality \(\le\aleph_2\).

We shall consider models \(\ma,\mb\) of cardinality \(\aleph_2\),
so we may as well assume they have \(\go_2\) as universe. For such
a model \(\ma\) and any ordinal \(\ga<\go_2\) we let \(\ma_{\ga}\)
denote the structure \(\ma\cap\ga\). Similarly \(\mb_{\ga}\). Let
us first recall the following basic fact from \cite{msv416}:

\begin{lemma}\cite{msv416}
\label{stationary} Suppose \(\ma\) and \(\mb\) are structures of
cardinality \(\aleph_2\). If \(\forall\) does not have a winning
strategy in \(\EFG{\goy}(\ma,\mb)\), then
\[S=\{\ga : \ma_{\ga}\cong\mb_{\ga}\}\]
is \(\goy\)-stationary.
\end{lemma}

This shows that to get determinacy of \(\EFG{\goy}(\ma,\mb)\) it
suffices to give a winning strategy of \(\exists\) under the
assumption that the above set \(S\) is \(\goy\)-stationary. In
\cite{msv416} an assumption \(I^*(\omega)\) was used. This
assumption says that the non-\(\goy\)-stationary ideal on
\(\omega_2\) has a \(\sigma\)-closed dense set. The rough idea was
that \(\exists\) uses the Pressing Down Lemma on \(S\) to
"normalize" his moves so that he always has an \(\goy\)-stationary
sets of possible continuations of the game. We use now the same
idea. The hypothesis \(I^*(\omega)\) is equiconsistent with a
measurable cardinal. Since we assume only the consistency of a
weakly compact cardinal, we have to work more.

Suppose \(\kappa\) is a weakly compact cardinal. Let \(\I\) denote
the
 \(\Pi^1_1\)-ideal on \(\kappa\),
i.e..
 the ideal
of subsets of \(\kappa\) generated by the sets
\(\{\alpha:(H(\alpha),\epsilon,A\cap
H(\alpha))\models\neg\phi\}\), where  \(A\subseteq H(\kappa)\) and
\(\phi\) is a \(\Pi^1_1\)-sentence such that
\((H(\kappa),\epsilon,A)\models\phi\). We collapse \(\kappa\) to
\(\omega_2\) and then force a cub to the complement of every set
\(S\subseteq S^2_1\) in \(\I\). In the resulting model the above
"normalization" strategy of \(\exists\) works even though the
non-\(\goy\)-stationary ideal on \(\omega_2\) may not have a
\(\sigma\)-closed dense set.

\begin{definition}
Let \(\F\) be a set of cardinality \(\kappa\) of regressive
functions \(\kappa\rightarrow\kappa\) and
\(S\subseteq\kappa\). The game ${\BM_{\omega_1}}({S,\F})$ has two
players called \(\emp\) and \(\non\). They alternately play
\(\omega_1\) rounds. During each round  \(\emp\) first chooses
\(f_i\in\F\). Then \(\non\)
chooses a subset \(S_i\) of \(\bigcap_{j<i}S_j\) (of
\(S\), if \(i=0\)) such that
it is unbounded in $\k$ and
\(f_{i}\) is constant on \(S_i\).
Player \(\non\) wins if he can play all \(\omega_1\) moves
following the rules.
\end{definition}

%



\begin{lemma}\label{wtyu}
Suppose \(S=\{\ga<\omega_2 : \ga\ne 0, \ma_{\ga}\cong\mb_{\ga}\}\) and
\(h_\alpha:\ma_\alpha\cong\mb_\alpha\) for \(\alpha\in S\).
Let \[\F=\{f_\alpha:\alpha\in S\}\cup
\{g_\alpha:\alpha\in S\},\] where 
\(f_\alpha:\omega_2\rightarrow\omega_2\) is the regressive
function mapping \(\xi\) (\(\ne 0\)) to \(h_\xi(\alpha)\) if \(\xi
>\alpha\), and to \(0\) otherwise, and \(g_\alpha\)
is the regressive
function mapping \(\xi\) (\(\ne 0\)) to \((h_\xi)^{-1}(\alpha)\) 
if \(\xi>\alpha\), and to \(0\) otherwise.
Suppose \(\non\) has a winning strategy in
\(\BM_{\omega_1}(S,\F)\). Then \(\exists\) has a winning strategy
in the game \(\EFG{\goy}^{\aleph_2}(\ma,\mb)\).
\end{lemma}

\proof We present the proof for \(\EFG{\goy}^{2}(\ma,\mb)\). The
case of \(\EFG{\goy}^{\aleph_2}(\ma,\mb)\) is similar. 
\(\H=\{h_\alpha:\alpha\in S\}\), where
\(h_\alpha:\ma_\alpha\cong\mb_\alpha\) for \(\alpha\in S\).
 Let
\(\tau\) be a winning strategy of \(\non\) in the game
\(\BM_{\omega_1}(S,\F)\). Suppose the sequence \(\la (x_i,y_i) :
i<\ga\ra\) has been played, where \(\alpha<\omega_1\), \(x_i\)
denotes a move of \(\forall\) and \(y_i\) a move of \(\exists\).
Suppose \(\forall\) plays next \(x_{\ga}\). 
During the game \(\exists\) also plays
\(\BM_{\omega_1}(S,\F)\). Let us denote his moves in
\(\BM_{\omega_1}(S,\F)\) by \(S_i\). Thus \(S_j\subseteq S_i\) for
\(i<j<\alpha\). The point of the sets \(S_i\) is that \(\exists\)
has taken care that for all \(i<\ga\) and \(j\in S_i\) we have
\(y_i=h_{j}(x_i)\) or \(x_i=h_{j}(y_i)\) depending on whether
\(x_i\in\ma\) or \(x_i\in\mb\). Let
$S'_{\ga}=\bigcap_{i<\alpha}S_i\setminus\alpha_i$. 
 The winning strategy $\tau$
gives an \(S_{\ga}\se S'_{\ga}\)
 and a \(y_{\ga}\) such that
\(f_i(x_{\ga})=y_{\ga}\) for all \(i\in S_{\ga}\),
if  \(x_{\ga}\in\ma\), and 
\(g_i(x_{\ga})=y_{\ga}\) for all \(i\in S_{\ga}\),
if  \(x_\alpha\in\mb\). This element
\(y_{\ga}\) is the next move of \(\exists\). Using this strategy
\(\exists\) cannot lose and hence wins. \qed



\begin{theorem}\label{13}
It is consistent relative to the consistency of a weakly compact
cardinal, that for every \(\goy\)-stationary
\(S\subseteq\omega_2\) and every set \(\F\) of cardinality
\(\aleph_2\) of regressive functions \(\omega_2\rightarrow\omega_2\),
$\non$ has a winning strategy in the game
$\BM_{\omega_1}(S,\F)$.
\end{theorem}

\proof 
We may assume GCH. Suppose \(\kappa\) is weakly 
compact.
Let \(\QQ\) be the Levy-collapse of \(\kappa\) to \(\aleph_2\). In
\(V^{\QQ}\) we define by induction a sequence \(\PP_\alpha\),
\(\alpha<\kappa^+\), of forcing notions. Let \((A_\alpha),
{\alpha<\kappa^+}\),
be a complete list of all sets in the $\Pi^1_1$-ideal \(\I\) on
$\kappa$ such that every element of $A_\alpha$ has uncountable
cofinality. If \(\alpha\) is limit of cofinality \(\le\omega_1\),
then \(\PP_\alpha\) is the inverse limit of all \(\PP_\beta\),
\(\beta<\alpha\). For other limit \(\alpha\), \(\PP_\alpha\) is
the direct limit of \(\PP_\beta\), \(\beta<\alpha\). At successor
stages we let \(\PP_{\alpha+1}=\PP_\alpha\star \R_\alpha\),
where \(\R_\alpha\) is defined as follows: \(q\in
\R_\alpha\) iff \(q\) is
a bounded closed sequence of elements of \(\kappa\) such that
\( q\cap {A_\alpha}=\emptyset\). \(\R_{\a}\) is ordered by the
end extension relation.
Thus each \(\PP_\alpha\) is countably closed.
Let \(\PP=\PP_{\kappa^+}\). Now \(\QQ\star\PP\) satisfies the
\(\kappa^+\)-chain condition. Note also that for all
$\a <\k^{+}$, $\QQ\star\PP_{\a}$ has power $\k$.
We prove that it is true in
\(V^{\QQ}\) that \(\PP_\alpha\) does not add new subsets of
\(\kappa\) of cardinality $\le\aleph_1$, hence \(\kappa\) remains \(\aleph_2\) also after
forcing with \(\PP\). It follows also that $\QQ\star\PP$ and
each $\QQ\star\PP_\alpha$ are countably closed.

We show now that in \(V^{\QQ\star\PP}\) the claim is true. Suppose
\(S\) and a set  \(\F=\{f_\alpha:\alpha<\kappa\}\),
of regressive functions \(\kappa\rightarrow\kappa\) are given in
\(V^{\QQ\star\PP}\) such that (in \(V^{\QQ\star\PP}\))
\(S\subseteq S^2_1\) is $\omega_1$-stationary. Suppose
\(\alpha<\kappa^+\) is such that \(\tilde{S},\tilde{\F}\) and
\(\tilde{f_i}\) are \(\QQ\star\PP_\alpha\)-names for \(S,\F\)
and \(f_i\), correspondingly.
Since \(S\) is  \(\omega_1\)-stationary in \(V^{\QQ\star\PP}\),
\(S\) is not in the ideal generated by \(\I\) in
\(V^{\QQ\star\PP_\alpha}\).
Suppose
\((p,q)\force
\tilde{S}\notin\I\). For a contradiction, suppose also that
$(p,q)$ forces that \(\non\) does not have a winning strategy
in the game \(BM_{\omega_1}(S,\F)\).

Let \((\B'_0,\in )\) be a sufficiently elementary substructure
of \((V,\in )\) such
that \(|\B'_0|=\kappa\), \({\B'}_0^{<\kappa}\se\B'_0\), \(\QQ\),
\(\PP_\alpha\), \(\alpha\), \(\kappa\), \(\tilde{\F},
\tilde{f_i}\), \(\tilde{S}\), are in \(\B'_0\), and
$\alpha\cup\kappa\subseteq\B'_0$. Let $\B_0$ be the transitive
collapse of $\B'_0$. Thus $\QQ ,\PP_\a ,\a ,\k\in\B_0$,
$A_j\in \B_0$ for $i\le\alpha$ and
$\tilde{f_i}\in\B_0$ for $i<\kappa$.
Let
$$T=\{\a <\k\vert\ \exists (p',q')\le (p,q)
((p',q')\force_{\QQ\star\PP_\alpha}\a\in\tilde{S})\}.$$
Clearly $T\in\B_0$ and $T\not\in\I$.
By weak compactness, there are a
transitive \(\B_1\) and an elementary embedding
\(j:\B_0\rightarrow\B_1\)
such that \(\kappa\) is the critical point of \(j\),
$\k\in j(T)$ and $\k\notin j(A_i)$ for
\(i\le\alpha\). So there is some $(p',q')\in j(\QQ\star\PP_\a )$
such that $(p',q')\le j((p,q))=(p,q)$ and
$(p',q')\force_{j(\QQ\star\PP_\a )}\k\in j(\tilde{S})$.
Note that
$\QQ ,\PP_\a\in\B_1$ and
$\tilde{f_i}\in\B_1$ for $i<\kappa$.

By (the proof of) Lemma 3 in \cite{mag}, there are
a $\QQ\star\PP_\a$-generic $G$ over $\B_1$ and a forcing notion
$\R\in\B_1[G]$ such that $(p,q)\in G$,
in $\B_1[G]$,
$\R$ is countably closed, for all $\R$-generic $K$ over
$\B_1[G]$, there is a canonical
\(j(\QQ\star\PP_\alpha )\)-generic $G_K$
over $\B_1$ such that $\B_1[G_K]=\B_1[G][K]$ and for some
$K$, $G_K$ is such that $(p',q')\in G'$.
Then
for every \(\QQ\star\PP_\alpha\)-name $\tilde{X}\in\B_0$,
there is a canonical $\R$-name $\tilde{Y}\in\B_1[G]$ such that for all
$\R$-generic $K$ over $\B_1[G]$, $j(\tilde{X})$ and $\tilde{Y}$
have the same interpretation in $\B_1[G][K]$.
We do not distinguish $j(\tilde{X})$
and $\tilde{Y}$.
With this notation,
there is $r\in\R$ which forces in $\B_1[G]$, that $\k\in j(\tilde{S})$.
Then there is some $(p^{*},q^{*})\le (p,q)$ in $G$ that
in $\B_1$ forces the existence of such $\R$ and $r$.
So we may assume that $G$ is generic over $V$ and our 
$V^{\QQ\star\PP_\a}$
is the same as $V[G]$.

We describe in $\B_1[G]$ a winning
strategy of \(\non\) in the game \(BM_{\omega_1}(S,\F)\).
This is a contradiction since all possible winning plays of
$\emp$ are in $\B_1[G]$ and being unbounded is absolute
in transitive models.
The strategy of $\non$
is to play on the side conditions
\(q^i\) in $\B_1[G]$ and sets \(S_i\in\B_0[G]\)
with \(\QQ\star\PP_\alpha\)-names
\(\tilde{S_i}\) in $\B_0$
such that
\begin{enumerate}
\item \(q^i\in\R\).
\item \(q^0\le r\).
\item \(i<k<\omega_1\) implies
\(q^k\le q^i\).
\item \(i<k<\omega_1\) implies
\(S_k\subseteq S_i\subseteq S\).
\item \(q^i\force_{\R} \kappa\in j(\tilde{S_i})\) in \(\B_1[G]\).
\end{enumerate}

Suppose \(\non\) has followed this strategy, forming
conditions \(q^i\)  and sets \(S_i\) for \(i<k\).
Let \(p=
\inf(\{q^i:i<k\})\). If we let \(S\) to be
\(\bigcap_{i<k}S_i\) and \(\tilde{S}\) a name for this, then
in $\B_1[G]$,
\[p\force_{\R} \kappa\in j(\tilde{S}).\]
%
%
%
Suppose
then \(\emp\) moves \(f_k\in\F\). 
Let $q^k\le p$ such that for some
$\delta<\kappa$ we have $q^k\force_{\R}
j(\tilde{f_k})(\kappa)=\delta$ in $\B_1[G]$ and let $S_k$ be
$\{\beta\in S:f_k(\beta)=\delta\}$ and $\tilde{S_k}$ a name for this.
Then $q^k\force_{\R}\kappa\in j(\tilde{S_k})$ in $\B_1[G]$.

Finally we have to prove that $\QQ\star\PP_\alpha$ does not add
new  subsets of \(\kappa\) of cardinality $\le\aleph_1$ 
over and above those added by $\QQ$. 
The proof
of this is, mutatis mutandis, like the proof of the Main fact
(page 761) in \cite{mag}.
Here we use the assumption $\kappa\notin
j(A_i)$ for $i\le\alpha$. Thus, if $C$ is a generic sequence
in the complement of $j(A_\beta)$ in $V^{j(\QQ\star\PP_\alpha)}$,
then we can continue it to a closed condition
$C\cup\{\kappa\}\in\R_{j(\beta)}$.
\qed

Results similar to Theorem~\ref{13} have been treated also in
\cite{SS87} and \cite{SS93}.

\begin{corollary}\label{i}
It is consistent relative to the consistency of a weakly compact
cardinal, that the game \(\EFG{\omega_1}(\ma,\mb)\) is determined
for all \(\ma\) and \(\mb\) of cardinalty \(\le\aleph_2\).
\end{corollary}



\section{Non-determinacy and structure theory}

In this section we prove Theorem~\ref{1}, which essentially
establishes, under cardinality assumptions concerning the
continuum, the existence of non-determined \ef games of length
\(\omega_1\) for models of {\em non-classifiable} theories. This
complements the observation, made in \cite{msv416}, that the \ef
game of length \(\omega_1\) is determined for models of {\em
classifiable} theories.

We start be proving Theorem~\ref{1} under assumption (iii), which
we consider the most intereting case. That is, we start with a
countable complete stable and unsuperstable first order theory and
show that, assuming $2^{\o}\le\o_{3}$, it has two models \(\ma\)
and \(\mb\) of cardinality \(\aleph_3\) for which
\(\EFG{\omega_1}(\ma,\mb)\) is non-determined. Actually, we
construct \(\ma\) and \(\mb\) so that \(\exists\) does not have a
winning strategy even in \(\EFG{\omega+\omega}^2(\ma,\mb)\) and
\(\forall\) does not have a winning strategy even in
\(\EFG{\omega_1}^{\omega_3}(\ma,\mb)\).

We then prove Theorem~\ref{1} under assumption (i), that is, we
now start with a countable complete unstable  first order theory
and show that, assuming $2^{\o}<2^{\o_{3}}$, it has two models
\(\ma\) and \(\mb\) of cardinality \(\aleph_3\) for which
\(\EFG{\omega_1}(\ma,\mb)\) is non-determined.

Theorem~\ref{1} under assumption (ii) can be dealt with in
 the same way
as under assumption (i). The section ends with some remarks on
possible improvements.

\subsection{The stable unsuperstable case}

We will prove Theorem~\ref{1}, case (iii), in a series of lemmas.
We assume $\o_{3}^{\o}=\o_{3}$ all the time. Let \(T\) be a
countable complete stable and unsuperstable first order theory. As
usual, we work inside a large saturated model $\M$ of $T$. We
start by fixing some notation. By a tree $I$ we mean a
lexicographically ordered downwards closed subtree of
$\theta^{<(\o +1)}$ for some linear order $\theta$,  that is,
$I=(I,\ll ,P_{\a},<,H)_{\a\le\o}\in K^{\o}_{tr}(\theta )$, see
\cite{ht} Definition 8.2 or \cite{Sh:c}. For a while, we fix a
tree $I\in K^{\o}_{tr}(\l )$, where $\l$ is some large enough
cardinal, so that $(I,\ll )$ is isomorphic to $\l^{<(\o +1)}$. As
in \cite{HS}, for $u,v\in \P_{\o}(I)$ (=finite subsets of $I$), we
define $r(u,v)$ to be the unique set $R$ which satisfies

\begin{description}
\item[(I)] $R\subseteq X_{u,v}=\{ H(\n ,\xi )\vert\ \n\in u,\ \xi\in v\}$,

\item[(II)] For all $\nu\in X_{u,v}-R$, there is $\nu'\in R$
such that $\nu\ll\nu'$,

\item[(III)] If $\n$ and $\xi$ are distinct elements of $R$, then
$\n\not\ll\xi$.
\end{description}

\noindent We write $u\le v$ if $r(u,v)=r(u,u)$. For more on these
definitions, see \cite{HS}. In \cite{HS}, it is shown that there
are models $\A$ and $\A_{u}$, $u\in\P_{\o}(I)$, and sequences
$a_{\n}$ from $\A_{\{\n\}}$, $\n\in I$, such that

\begin{description}
\item[(i)] $\A =\bigcup_{u\in\P_{\o}(I)}\A_{u}\models T$,

\item[(ii)] if $u\le v$, then $\A_{u}\subseteq\A_{v}$,

\item[(iii)] for all $u,v\in\P_{\o}(I)$,
$\A_{u}\da_{\A_{r(u,v)}}\A_{v}$,

\item[(iv)] for all $u\in\P_{\o}(I)$, $\vert\A_{u}\vert\le\o_{3}$,

\item[(v)] if $P_{\o}(\n )$ holds and $\xi\ll\n$ is an immediate
successor of $\xi'$, then $$a_{\n}\nda_{\A_{\{\xi'\}}}a_{\xi}.$$
\end{description}

\noindent These models are exactly what we want except that they
are too large, we want the models $\A_{u}$, $u\in\P_{\o}(I)$, to
be countable. In order to get this, we use the
Ehrenfeucht-Mostowski construction.

We extend the signature $L$ of $T$ to $L_{*}$ by adding $\o_{3}$
new function symbols, some of which will be interpreted in $\M$ so
that they provide Skolem-functions for the $L$-formulas. In
addition we interpret the functions so that if we write
$SH_{*}(u)$ for the $L_{*}$-Skolem-hull of $\{ a_{\n}\vert\
\{\n\}\le u\}$ then
\begin{description}
\item[(vi)] for all $u\in\P_{\o}(I)$, $SH_{*}(u)=\A_{u}$.
\end{description}
By the usual argument (using \cite[Appendix Theorem 2.6]{Sh:c} and
compactness) we can interpret the new function symbols so that
$\M$ remains sufficiently saturated and the following holds
\begin{description}
\item[(vii)] if $U$ is a downwards closed subtree of $I$ and
$f$ is an automorphism of $U$, then there is an
$L_{*}$-automorphism $g$ of $\cup_{u\in\P_{\o}(U)}\A_{u}$ such
that for all $\n\in U$, $g(a_{\n})=a_{f(\n )}$.
\end{description}
\noindent Finally, it is easy to see that we can choose countable
$L_{1}\subseteq L_{*}$ so that $L\subseteq L_{1}$, $L_{1}$
contains the Skolem-functions for the $L$-formulas and if we write
$SH_{1}(u)$ for the $L_{1}$-Skolem-hull of $\{ a_{\n}\vert\
\{\n\}\le u\}$ then
\begin{description}
\item[(viii)] for all $u,v\in\P_{\o}(I)$, $SH_{1}(u)\da_{SH_{1}(v)}SH_{*}(v)$.
\end{description}
\noindent So we have proved the following lemma (for the notion
$\Phi$ proper for $K^{\o}_{tr}$ and the Ehrenfeucht-Mostowski
models $EM^{1}(J,\Phi )$, see \cite{ht} Definition 8.1 or
\cite{Sh:c}).

\begin{lemma}\label{2.3}
There are countable $L_{1}\supseteq L$ and $\Phi$ proper for
$K^{\o}_{tr}$ such that the following holds:
\begin{description}
\item[(a)] For all $J\in K^{\o}_{tr}$ there are an $L_{1}$-model
$EM^{1}(J,\Phi )\models T$ and sequences $a_{\n}\in EM^{1}(J,\Phi
)$, $\n\in J$, such that $EM^{1}(J,\Phi )$ is the
$L_{1}$-Skolem-hull of $\{ a_{\n}\vert\ \n\in J\}$ (i.e. $\{
a_{\n}\vert\ \n\in J\}$ is the skeleton of $EM^{1}(J,\Phi )$ and
as before for $u\subseteq J$, $SH_{1}(u)$ denotes the
$L_{1}$-Skolem hull of $\{ a_{\n}\vert\ \{\n\}\le u\}$).

\item[(b)] If $U$ is a downwards closed subtree of $J$ and
$f$ is an automorphism of $U$, then there is an
$L_{1}$-automorphism $g$ of $SH_{1}(U)$ such that for all $\n\in
U$, $g(a_{\n})=a_{f(\n )}$.

\item[(c)] Assume $(\n_{i})_{i<\o}$ is a strictly $\ll$-increasing
sequence of elements of $J$, $\n_{i+1}$ is an immediate successor
of $\n_{i}$ and $\n_{0}$ is the root. Then $(\n_{i})_{i<\o}$ has
an upper bound in $J$ iff there is a sequence $a\in EM(J,\Phi )$
such that for all $i<\o$, $a\nda_{SH_{1}(\{\n_{i}\}
)}a_{\n_{i+1}}$. $\eop$

\end{description}
\end{lemma}

We will write $EM(J,\Phi )$ for $EM^{1}(J,\Phi )\raj L$.

Our next goal is to define the skeletons for the models $\A$ and
$\B$ in the theorem. For this we use the weak box from
\cite{msv416}. By $S^{n}_{m}$ we denote the set $\{\a<\o_{n}\vert\
cf(\a )=\o_{m}\}$.

\begin{theorem}\label{2.4}
(\cite[Lemma 16]{msv416}) There are sets $S$, $U$ and $C_{\a}$,
$\a\in S$, such that the following holds:

(a) $S\subseteq S^{3}_{0}\cup S^{3}_{1}$ and $S\cap S^{3}_{1}$ is
stationary,

(b) $U\subseteq S^{3}_{0}$ is stationary and $S\cap U=\empty$,

(c) for all $\a\in S$, $C_{\a}\subseteq\a\cap S$ is closed in $\a$
and of order-type $\le\o_{1}$,

(d) for all $\a\in S$, if $\b\in C_{\a}$, then
$C_{\b}=C_{\a}\cap\b$,

(e) for all $\a\in S\cap S^{3}_{1}$, $C_{\a}$ is unbounded in
$\a$. $\eop$

\end{theorem}

We will construct trees $I_{\a}$ and $J_{\a}$, $\a <\o_{3}$, so
that the following holds:

\begin{description}
\item[(1)] if $\a <\b$ then $I_{\a}$ is a submodel of $I_{\b}$ and
$J_{\a}$ is a submodel of $J_{\b}$; now for $\n\in I_{\a}$, we
will write $rk(\n )$ for the least $\b$ such that $\n\in I_{\b}$
and similarly for $\n\in J_{\a}$,

\item[(2)] for all $\a\in S$, there is an isomorphism
$G_{\a}:I_{\a}\rightarrow J_{\a}$,

\item[(3)] if $\a\in C_{\b}$, then $G_{\a}\subseteq G_{\b}$,

\item[(4)] for all $\a\le\b$ and $\n\in I_{\a}$, if
$P_{\o}(\n )$ does not hold, then there is an immediate successor
$\xi$ of $\n$ such that $\xi\in I_{\b +1}-I_{\b}$

\item[(5)] if $(\n_{i})_{i<\o}$ is an increasing sequence of
elements of $I_{\a}$ (for some $\a$) and the sequence has an upper
bound $\xi$ in $I_{\a}$, then $rk(\xi )=sup_{i<\o}rk(\n_{i})$ and
similarly for sequences from $J_{\a}$,

\item[(6)] if $(\n_{i})_{i<\o}$ is an increasing sequence of
elements of $I_{\a}$, $(rk(\n_{i}))_{i<\o}$ is not eventually
constant and the sequence has an upper bound $\xi$ in $I_{\a}$,
then $rk(\xi )$ $(=sup_{i<\o}rk(\n_{i}))$ $\in U$; in $J_{\a}$
such sequences never have an upper bound,

\item[(7)] $\vert I_{\a}\vert\le\o_{3}$ and $\vert J_{\a}\vert\le\o_{3}$,

\item[(8)] $I_{\a},J_{\a}\subseteq H_{\o}(\o_{3})$, where $H_{\o}(\o_{3})$
is the least set $H$ such that $\o_{3}\subseteq H$ and if
$E\subseteq H$ is of power $\le\o$, then $E\in H$.
\end{description}

\noindent It is easy to see that such trees can be constructed by
induction on $\a$. However, in order to get what we want we need
to do a bit more work when we define $I_{\a}$ and $J_{\a}$ in the
case $\a\in U$. In order to decide, which branches like the one in
(6) above, we want to have an upper bound, we use a guessing
machine from \cite{Sh2} called black box, which we formulate so
that it fits exactly to our purposes.

\begin{theorem}\label{2.5}
 (\cite{Sh2}) ($\o_{3}^{\o}=\o_{3}$.)
There are $(\ol M^{\a},\n^{\a})$, $\a <\o_{3}$, such that

\begin{description}
\item[(i)] $\ol M^{\a}=(M^{\a}_{i})_{i<\o}$ is an increasing elementary
chain of elementary submodels of some $(H_{\o}(\o_{3}),A,B,\s )$,
such that $A,B\subseteq H_{\o}(\o_{3})$ and $\s$ is a strategy of
$\exists$ in $EF^{2}_{\o}(A,B)$ ($A$ and $B$ can be viewed as models of
the empty signature),

\item[(ii)] $M^{\a}_{i}=(M^{\a}_{i},A^{\a}_{i},B^{\a}_{i},\s^{\a}_{i})\in
H_{\o}(\o_{3})$,

\item[(iii)] $\n^{\a}$ is an increasing function from $\o$ to $\o_{3}$,
$\M^{\a}_{i}\in H_{\o}(\n^{\a}(i+1))$ and $sup_{i<\o}\n^{\a}(i)\in
U$,

\item[(iv)] $(\n^{\a}(j))_{j\le i},(M^{\a}_{j})_{j\le i}\in M^{\a}_{i+1}$,

\item[(v)] if $\a\ne\b$, then $\n^{\a}\ne\n^{\b}$,

\item[(vi)] player I does not have a winning strategy for the following game:
The length of the game is $\o$. At each move $i<\o$, first I
chooses $M_{i}$ and then II chooses $\a_{i}<\o_{3}$. I must play
so that in the end (i), (ii) and (iv) above are satisfied. I wins
if he has played according to the rules and there is no $\a
<\o_{3}$ such that $((M_{i})_{i<\o},(\a_{i})_{i<\o})=(\ol
M^{\a},\n^{\a})$.
\end{description}
\end{theorem}

First we uniformize (partially) the Ehrenfeuct-Mostowski
construction: We assume that for all $I,I'\in K^{\o}_{tr}$, if $I$
is a substructure of $I'$ and $I'\subseteq H_{\o}(\o_{3})$, then
there is a unique model $EM^{1}(I,\Phi )$, it is a substructure of
$EM^{1}(I',\Phi )$ and $EM^{1}(I',\Phi )\subseteq H_{\o}(\o_{3})$.

So let $\a\in U$ and assume that $I_{\b}$ and $J_{\b}$ are defined
for all $\b <\a$. Write $I^{*}_{\a}=\cup_{\b <\a}I_{\b}$ and
$J^{*}_{\a}=\cup_{\b <\a}I_{\b}$. For $\g <\o_{3}$, we write
$M^{\g}$ for $\cup_{i<\o}M^{\g}_{i}$ and $A^{\g}$, $\B^{\g}$ and
$\s^{\g}$ are defined similarly. Let $W^{\a}$ be the set of all
$\g <\o_{3}$ such that
\begin{description}
\item[(a)] $A^{\g}=EM(I^{*}_{\a}\cap M^{\g},\Phi )$ and
$B^{\g}=EM(J^{*}_{\a}\cap M^{\g},\Phi )$,

\item[(b)] $sup_{i<\o}\n^{\g}(i)=\a$,

\item[(c)] there are $\xi^{\g}_{i}\in I^{*}_{\a}\cap M^{\g}$,
$i<\o$, such that $\xi^{\g}_{0}$ is the root of $I^{*}_{\a}$,
$\xi^{\g}_{i+1}$ is an immediate successor of $\xi^{\g}_{i}$ and
$\xi^{\g}_{i}\in I_{\n^{\g}(i)+1}-I_{\n^{\g}(i)}$.
\end{description}

\noindent Notice that by Theorem~\ref{2.5} (v), if $\g\ne\d$, then
$(\xi^{\g}_{i})_{i<\o}\ne (\xi^{\d}_{i})_{i<\o}$. Let
$C^{\g}_{i}=SH(\{\xi^{\g}_{i}\} )$. Then we can find a partial
function $g^{\g}:A^{\g}\rightarrow\B^{\g}$ such that
\begin{description}
\item[(d)] $dom(g^{\g})=\cup_{i<\o}C^{\g}_{i}$,

\item[(e)] $g^{\g}$ is a result of a play of $EF^{2}_{\o}(A^{\g},B^{\g})$
in which $\exists$ has used $\s^{\g}$.
\end{description}
\noindent We let $W^{\a}_{J}$ be the set of those $\g\in W^{\a}$
such that
\begin{description}
\item[(f)] $g^{\g}$ is a partial isomorphism from
$EM(I^{*}_{\a}\cap M^{\g},\Phi )$ to $EM(J^{\g}_{\a}\cap
M^{\g},\Phi )$,

\item[(g)] there is $J$ such that if we let $J_{\a}=J$, then
(1),(5)-(8) above are satisfied and there is a sequence $a\in
EM(J,\Phi )$ such that for all $i<\o$,
$a\nda_{g^{\g}(C^{\g}_{i})}g^{\g}(a_{\xi^{\g}_{i+1}})$.
\end{description}
\noindent We let $W^{\a}_{I}$ be the set of all $\g\in
W^{\a}-W^{\a}_{J}$ such that $g^{\g}$ satisfies (f) above.

Now we can define $I_{\a}$ and $J_{\a}$. First we choose $I_{\a}$
so that it consists of all $\n\in I^{*}_{\a}$ together with the
supremums for the branches $(\xi^{\g}_{i})_{i<\o}$, $\g\in
W^{\a}_{I}$. $J_{\a}$ is chosen so that it satisfies (g) for all
$\g\in W^{\a}_{J}$ (and so especially (1),(5)-(8)).

Then we let $I=\cup_{\a <\o_{3}}I_{\a}$, $J=\cup_{\a
<\o_{3}}J_{\a}$, $\A =EM(I,\Phi )$ and $\B =EM(J,\Phi )$. Clearly
$\A$ and $\B$ can be chosen so that $\A ,\B\subseteq
H_{\o}(\o_{3})$.

\begin{lemma}\label{2.6}
$\forall$ does not have a winning strategy for $EF^{\o_{3}}_{\o_{1}}(\A
,\B )$.
\end{lemma}

\proof For this it is enough to show that $A$ does not have a
winning strategy for $EF^{\o_{3}}_{\o_{1}}(I,J)$, which is clear
by (2) and (3) above and Theorem~\ref{2.4}. $\eop$

\begin{lemma}\label{2.7}
$\exists$ does not have a winning strategy for $EF^{2}_{\o +\o}(\A ,\B
)$.
\end{lemma}

\proof For a contradiction, assume $\s$ is a winning strategy of
$\exists$ for the game $EF^{2}_{\o +\o}(\A ,\B )$. We play a round of
the game defined in Theorem~\ref{2.5} (vi). We let player I play
so that he follows the rules and
\begin{description}
\item[(i)] for all $i<\o$, $M_{i}\prec (H_{\o}(\o_{3}),\A ,\B ,\s\raj\o )$,

\item[(ii)] for all $\d ,\d'\in M_{i}$, if $\d\le\d'$, $\n\in I_{\d}\cap M_{i}$
and $P_{\o}(\n )$ does not hold, then there is $\xi\in
(I_{\d'+1}-I_{\d'})\cap M_{i+1}$ such that $\xi$ is an immediate
successor of $\n$,

\item[(iii)] the Skolem-hulls of $\{ a_{\n}\vert\ \n\in I\cap M_{i}\}$
and $\{ a_{\n}\vert\ \n\in J\cap M_{i}\}$ are subsets of
$M_{i+1}$,

\item[(iv)] $\A\cap M_{i}$ is a subset of the Skolem hull of
$\{ a_{\n}\vert\ \n\in I\cap M_{i+1}\}$ and $\B\cap M_{i}$ is a
subset of the Skolem hull of $\{ a_{\n}\vert\ \n\in J\cap
M_{i+1}\}$,

\item[(v)] $\bigcup\{ rk(\n )\vert\ \n\in I\cap M_{i}\}\cup
\bigcup\{ rk(\n )\vert\ \n\in J\cap M_{i}\}\in M_{i+1}$.
\end{description}
\noindent By Theorem~\ref{2.5} (vi), the round can be played so
that $\forall$ loses. Let $\a_{i}$, $i<\o$, be the choices $\exists$ made and
$\g$ such that $((M_{i})_{i<\o},(\a_{i})_{i<\o})=(\ol
M^{\g},\n^{\g})$. Finally, let $\a =\cup_{i<\o}\a_{i}$ ($\in U$).

Now it is easy to see that $\g\in W^{\a}$, in fact $\g\in
W^{\a}_{I}$ or $\g\in W^{\a}_{J}$ (otherwise we have demonstrated
that $\s$ is not a winning strategy). In the first case, there is
a sequence $a\in\A$ such that for all $i<\o$,
$a\nda_{C^{\g}_{i}}a_{\xi^{\g}_{i+1}}$ but in $\B$ there is no
sequence $b$ such that for all $i<\o$,
$b\nda_{g^{\g}(C^{\g}_{i})}g^{\g}(\xi^{\g}_{i+1})$, a
contradiction. In the latter case, there is a sequence $b\in\B$
such that for all $i<\o$,
$b\nda_{g^{\g}(C^{\g}_{i})}g^{\g}(\xi^{\g}_{i+1})$ but by (the
construction,) Lemma 2.3 (c) and Theorem~\ref{2.5}  (v), there is
no sequence $a\in\A$ such that for all $i<\o$,
$a\nda_{C^{\g}_{i}}a_{\xi^{\g}_{i+1}}$, a contradiction. $\eop$

Now Lemmas 2.6 and 2.7 imply Theorem~\ref{1}  (iii).

\subsection{The unstable case}

We will prove Theorem~\ref{1}, case (i), again in a series of
lemmas. We assume $\o_{3}^{\o}<2^{\o_{3}}$. Let \(T\) be a
countable complete unstable first order theory. Let $L$ be the
signature of $T$.

\begin{theorem}\label{2.8}
(\cite{Sh:c}) Assume $T$ is a countable unstable theory in the
signarute $L$. There are a countable signature $L_{1}\supseteq L$,
a complete Skolem theory $T_{1}\supseteq T$ in the signature
$L_{1}$, a first-order $L$-formula $\phi (x,y)$ and $\Phi$ proper
for $(\o ,T_{1})$ (see [Sh1] Definition VII 2.6) such that for
every linear order $I$ there is an Ehrenfeucht-Mostowski model
$EM^{1}(I,\Phi )$ of $T_{1}$ with a skeleton $\{ a_{\n}\vert\
\n\in I\}$ such that $$EM^{1}(I,\Phi )\models\phi
(a_{\n},a_{\xi})\ \hbox{\sl iff} \ I\models\n <\xi .$$
\end{theorem}

We write $EM(I,\Phi )$ for $EM^{1}(I,\Phi )\raj L$. Notice that by
using the terminology from \cite[Definition III 3.1]{Sh2}, $\{
a_{\n}\vert\ \n\in I\}$ is weakly $(\o ,\phi)$-skeleton-like in
$EM(I,\Phi )$.

In order to use Theorem~\ref{2.8}, linear orders are needed. If
$A$ is a linear ordering, $x\in A$ and $B\subseteq A$, then by
$x<B$ we mean that for every $y\in B$, $x<y$, $x>B$ and $C>B$,
$C\subseteq A$ are defined similarly. By $A^{*}$ we mean the
inverse of $A$. Again let $S$, $U$ and $C_{\a}$, $\a\in S$, be as
in \cite[Lemma 16]{msv416}, i.e. Theorem~\ref{2.4} above, with the
exception that $0\in S$ and for all $\a\in S-\{ 0\}$, $0\in
C_{\a}$. By induction on $i<\o_{3}$, we will define linear orders
$A^{i}_{\a}$ and $B^{i}_{\a}$, $\a <\o_{3}$, and for $i\in S$,
isomorphisms $$G_{i}:\Sigma_{\b <i+2}A^{i}_{\b}\rightarrow
\Sigma_{\b <i+2}B^{i}_{\b} .$$ We write $A^{i}(\b ,\a )$ for
$\Sigma_{\b\le \g <\a}A^{i}_{\g}$ and similarly $B^{i}(\b ,\a )$.
We will do the construction so that the following holds:
\begin{description}
\item[(1)] $A^{0}_{\a}\cong\o^{*}$ for all $\a <\o_{3}$ and
if $\a\not\in U$, then $B^{0}_{\a}\cong\o^{*}$ and otherwise
$B^{0}_{\a}\cong (\o_{1})^{*}$,

\item[(2)] If $i<j$, then $A^{i}_{\a}\subseteq A^{j}_{\a}$ and
$B^{i}_{\a}\subseteq B^{j}_{\a}$ and otherwise the sets are
distinct and if $j\in C_{i}$, then $G_{j}\subseteq G_{i}$,

\item[(3)] if $cf(\a )=\o$, then
$A^{0}_{\a}$ is coinitial in $A^{i}_{\a}$ and similarly for $B$.
\end{description}

\noindent We will do this by induction on $i$. However, in order
to be able to show that (3) holds in each step, we need additional
machinery.

Let $C\in\{ A,B\}$. We say that $(I,J)$ is a $(C,i,\b )$-cut if
$I$ is an initial segment of $C^{i}_{\b}$ and $J=C^{i}_{\b}-I$. We
say that the cut is basic if $I=\empty$. We define a notion of
forbidden cut by induction on $i$ as follows (we should talk about
$i$-forbidden cuts, but $i$ is always clear from the context):

\begin{description}
\item[(a)] for all limit $\b$, the basic $(C,0,\b )$-cut is forbidden,

\item[(b)] if $(I,J)$ is a $(C,i,\b )$-cut, $j<i$ and
$(C^{j}_{\b}\cap I,C^{j}_{\b}\cap J)$ is forbidden, then $(I,J)$
is forbidden,

\item[(c)] if $(I,J)$ is a forbidden $(A,i,\b )$-cut,
$I^{*}=I\cup\bigcup_{\g<\b}A^{i}_{\g}$ and $G_{i}(I^{*})$ is not
bounded by any $x\in\cup_{\g<\d}B^{i}_{\g}$ but some $y\in
B^{i}_{\d}$ bounds it, then $(G_{i}(I^{*})\cap
B^{i}_{\d},B^{i}_{\d}-G_{i}(I^{*}))$ is forbidden and similarly
for $A$ and $B$ reversed (and $G_{i}$ replaced by $(G_{i})^{-1}$).
\end{description}

\noindent Now we can state the additional properties we want our
construction have. Let $E\in\{ A,B\}$, $i,\b <\o_{3}$ and $(I,J)$
be a $(E,i,\b )$-cut.
\begin{description}
\item[(4)] If $(I,J)$ is forbidden,
then there is no $j<\o_{3}$ and $x\in E^{j}_{\b}$ such that
$I<x<J$.

\item[(5)] Assume $(I,J)$ is forbidden and
$j\in S$ is such that $j<i$ and either $E^{j}_{\b}\cap I$ is
cofinal in $I$ or $E^{j}_{\b}\cap J$ is coinitial in J (we say
that $\empty$ is both cofinal and coinitial in $\empty$). Then
$(E^{j}_{\b}\cap I,E^{j}_{\b}\cap J)$ is forbidden.

\item[(6)] If $\b$ is successor, then $E^{0}_{\b}$ is coinitial in $E^{i}_{\b}$.
\end{description}

\begin{lemma}\label{2.9}
Let $E\in\{ A,B\}$.
\begin{description}
\item[(i)] For all $i,\b <\o_{3}$,
if (5) holds upto the stage $i$, then $(E^{i}_{\b},\empty )$ is
not forbidden and neither is $(\empty ,E^{i}_{\b})$, if $\b$ is
successor.

\item[(ii)] For limit $\b$, every basic $(E,i,\b )$-cut is forbidden.

\item[(iii)] The property (4) implies the property (3).

\item[(iv)] If $i+1 <\b$ and $(I,J)$ is a forbidden
$(E,i,\b )$-cut, then it is basic (and $\b$ is limit).
\end{description}
\end{lemma}

\proof Immediate. $\eop$

Now we are ready to do the construction: For $i=0$, the linear
orders are defined by (1) and we let $G_{0}$ be the only possible
one. Clearly (1)-(6) hold. If $i\not\in S$ or $sup\ C_{i}=i$, then
we let $A^{i}_{\a}=\cup_{j<i}A^{j}_{\a}$,
$B^{i}_{\a}=\cup_{j<i}B^{j}_{\a}$ and if $i\in S$ (and $sup\
C_{i}=i$), then $G_{i}=G\cup\bigcup_{j\in C_{i}}G_{j}$, where $G$
is the obvious isomorphism from $A^{i}(i,i+2)$ to $B^{i}(i,i+2)$
(both are isomorphic to $\o^{*}+\o^{*}$). Now (1), (2), (4) and
(6) hold trivially. By Lemma 2.9 (iii), (3) holds. For (5), assume
that $C\in\{ A,B\}$ and $(I,J)$ is a forbidden $(C,i,\b )$-cut.
Now the reason why $(I,J)$ is forbidden is (b) in the definition
of forbidden cut (if $i\not\in S$, then this is trivial and
otherwise by the definition of $G_{i}$, (c) does not give
forbidden cuts not forbidden by (b)). But then (5) follows
immediately from the induction assumption.

We are left with the case $i\in S$ and $j=sup\ C_{i}<i$. Notice
that now $j\in C_{i}$. Let $\a <j+2$ and $A\ne\empty$ be an
initial segment of $A^{j}_{\a}$. Let $A^{+}=A\cup\bigcup_{\g
<\a}A^{j}_{\g}$. Then there is the least $\b <j+2$ such that
$B^{+}=G_{j}(A^{+})\cap (\cup_{\g\le\b}B^{j}_{\b})=G_{j}(A^{+})$.
Let $A'=(A^{j}_{\a}\cup A^{j}_{\a +1})-A$, $B=G_{j}(A)\cap
B^{j}_{\b}$ and $B'=(B^{j}_{\b}\cup B^{j}_{\b +1})-B$. Assume that
at least one of $C'=\{ x\in\cup_{k<i}(A^{k}_{\a}\cup A^{k}_{\a
+1})\vert\ A<x<A'\}$ and $D'=\{ x\in\cup_{k<i}(B^{k}_{\b}\cup
B^{k}_{\b +1})\vert\ B<x<B'\}$ is non-empty. Then by the induction
assumption, $B\ne\empty$. Let $C$ be a copy of $C'$ and $D$ a copy
of $D'$. Then we define $A^{i}_{\a}$ so that it contains
$\cup_{k<i}A^{k}_{\a}$ and in each cut like above we add $D$ so
that $A<C'<D<A'$ and $B^{i}_{\b}$ is defined similarly but now
$B<C<D'<B'$ (this is possible by (6) in the induction assumption).
Then, by (4) in the induction assumption, we can find an
isomorphism $G'_{i}:\cup_{\a <j+2}A^{i}_{\a}\rightarrow\cup_{\a
<j+2}B^{i}_{\a}$. Notice that by (5) in the induction assumption,
for all $\d <i$, the $(A,\d ,\a )$-cut $(A^{\d}_{\a}-A(\d ),A(\d
))$ and $(B,\d ,\b )$-cut $(B(\d ),B^{\d}_{\b}-B(\d ))$ are not
forbidden, where $A(\d )=\{ x\in A^{\d}_{\a}\vert\ x>C'\}$ and
$B(\d )=\{ x\in B^{\d}_{\b}\vert\ x<D'\}$. So we have not violated
the property (4).

For all $\a >j+1$, we let $A^{i}_{\a}=\cup_{k<i}A^{k}_{\a}$ and
$B^{i}_{\a}$ is defined similarly. However we will still make
changes to $B^{i}_{j+1}$ and $A^{i}_{i+1}$! Let $A$ be a copy of
$B^{i}(j+3,i+2)$ and $B$ be a copy of $A^{i}(j+2,i+1)$.
Furthermore, extend $A^{i}_{i+1}$ so that there is an isomorphism
$g:A^{i}_{i+1}\rightarrow B^{i}_{j+2}$ such that
$g(A^{0}_{i+1})=B^{0}_{j+2}$ (this is not a problem since
$A^{0}_{i+1}=\cup_{k<i}A^{k}_{i+1}\cong\o^{*}\cong B^{0}_{j+2}$
and by Lemma 2.9 (iv), the sets $A^{k}_{i+1}$, $k<i$, do not
contain forbidden $(A,k,i+1)$-cuts; so we do not violate (4)).
Then we add $A$ to (the extended) $A^{i}_{i+1}$ as an end segment
and $B$ to $B^{i}_{j+1}$ as an end segment. By Lemma 2.9 (i), this
does not violated (4). Now it is easy to extend $G'_{i}$ to
$G_{i}$ so that $G_{i}(A^{i}(j+2,i+1))=B$,
$G_{i}(A^{i}_{i+1}-A)=B^{i}_{j+2}$ and $G_{i}(A)=B^{i}(j+3,i+2)$.

Now (1), (2) and (6) hold trivially, (4) is already shown to hold
and by Lemma 2.9 (iii), (4) implies (3). So we are left to show
that

\begin{lemma}\label{2.10}
(5) holds.

\end{lemma}

\proof Assume $(I,J)$ is a forbidden $(E,i,\b )$-cut, $E\in\{
A,B\}$, and $\d\in S$ is such that $\d <i$ and $E^{\d}_{\b}\cap J$
is coinitial in J, the other case is similar. If $\b\ge j+1$ and
both $J\cap(\cup_{k<i}A^{k}_{j+1})$ and
$J\cap(\cup_{k<i}B^{k}_{j+1})$ are empty, then the claim follows
easily from Lemma 2.9 and the induction assumption. So we assume
that this is not the case. If $(I,J)$ is forbidden because of (b)
in the definition of forbidden cut, the claim follows from the
induction assumption. So we assume that $E=B$ and there is a
forbidden $(A,i,\g)$-cut $(C,D)$ such that $(I,J)$ is forbidden by
(c) applied to $(C,D)$ (the case $A$ and $B$ reversed is
symmetric). Since $(I,J)$ is not forbidded by (b) in the
definition of forbidded cut, $(C,D)$ must be forbidden because of
it, i.e. for some $\a <i$, $(A^{\a}_{\g}\cap C,A^{\a}_{\g}\cap D)$
is a forbidden $(A,\a ,\g )$-cut.

If
\begin{description}
\item[(\(\star\))] For no
$y\in A^{\a}_{\g}\cap D$, $y<A^{j}_{\g}\cap D$,
\end{description}

\noindent then by the induction assumption, $(B^{j}_{\b}\cap
I,B^{j}_{\b}\cap J)$ is a forbidden $(B,j,\b )$-cut and the claim
follows from the definition of forbidden cut if $\d\ge j$ and from
(5) in the induction assumption if $\d <j$. So we assume that
(\(\star\)) fails. Let $y$ be the bound. Then $\empty\ne D'=\{
z\in A^{i}_{\g}\cap D\vert\ z\le y\} \subseteq
A^{i}_{\g}-dom(G_{j})$. So by the construction,
$G_{i}(D')\subseteq J-B^{\d}_{\b}$ and for all $x\in B^{\d}_{\b}$,
ether $x<G_{i}(D')$ or $x>G_{i}(D')$. By the choice of the cut
$(C,D)$, there can not be $x\in J\cap B^{\d}_{\b}$ such that
$x<G_{i}(D')$. But then $G_{i}(D')<J\cap B^{\d}_{\b}$, which
contradicts the assumption that $J\cap B^{\d}_{\b}$ is coinitial
in $J$. $\eop$

Let $A=\Sigma_{\a<\o_{3}}\cup_{i<\o_{3}}A^{i}_{\a}$ and
$B=\Sigma_{\a<\o_{3}}\cup_{i<\o_{3}}B^{i}_{\a}$. Notice that by
(1) and (3), $inv^{1}_{\o}(A)$ differs from $inv^{1}_{\o}(B)$ in a
stationary set which consists of ordinals of cofinality $\o$ (for
the definition of $inv^{n}_{\o}$, see \cite[Definition III
3.4]{Sh2}. Let $S_{\a}\subseteq S^{3}_{0}$, $i<2^{\o_{3}}$, be
stationary sets such that for $\a <\b <2^{\o_{3}}$,
$S_{\a}\triangle S_{\b}$ is stationary and define
$\Psi_{\a}=\Sigma_{\a <\o_{3}}\tau_{\a}$, where $\tau_{\a}=A^{*}$
if $\a\not\in S_{\a}$ and otherwise $\tau_{\a}=B^{*}$. Notice that
for $\a\ne\b$, $inv^{2}_{\o}(\Psi_{\a})$ differs from
$inv^{2}_{\o}(\Psi_{\b})$ in a stationary set which consists of
ordinals of cofinality $\o$.

Finally, let $\A_{\a}=EM((\Psi_{\a})^{*}\cdot\o_{1},\Phi)$.

\begin{lemma}\label{2.11}
For all $\a ,\b <2^{\o_{3}}$, A does not have a winning strategy
for $EF^{\o_{3}}_{\o_{1}}(\A_{\a},\A_{\b})$.

\end{lemma}

\proof For this, it is enough to show that A does not have a
winning strategy for
$EF^{\o_{3}}_{\o_{1}}((\Psi_{\a})^{*}\cdot\o_{1},(\Psi_{\b})^{*}\cdot\o_{1})$,
which follows easily from (2) in the construction of $A$ and $B$
and Theorem ~\ref{2.4} (see e.g. \cite[Claim 3 in the proof of
Theorem 17]{msv416}). $\eop$

\begin{lemma}\label{2.12}
There are $\a <\b <2^{\o_{3}}$ such that E does not have a winning
strategy for $EF^{2}_{\o_{1}}(\A_{\a},\A_{\b})$.
\end{lemma}

\proof By using the usual forcing notion, we collapse $\o_{3}$ to
an ordinal of power $\o_{1}$. Since this forcing notion does not
kill those stationary subsets of $\o_{3}$ which consist of
ordinals of cofinality $\o$ and cofinalities $\le\o_{1}$ are
preserved, in the extension, $inv^{2}_{\o}(\Psi_{\a})\ne
inv^{2}_{\o}(\Psi_{\b})$ for all $\a\ne\b$. Clearly, the skeletons
of the models $\A_{\a}$, remain weakly $(\o ,\phi )$-skeleton-like
in $\A_{\a}$. So by (the proof of) \cite[Lemma III 3.15 (1)]{Sh2},
$inv^{2}_{\o}(\Psi_{\a})\in INV^{2}_{\o}(\A_{\a},\phi )$ in the
extension (for the definition of $INV^{n}_{\o}$, see
\cite[Definition III 3.11]{Sh2} and notice that $\A\cong\B$
implies $INV^{2}_{\o}(\A ,\phi )=INV^{2}_{\o}(\B ,\phi )$). Also
by [Sh2] Lemma III 3.13 (1), $\vert INV^{2}_{\o}(\A_{\a},\phi
)\vert =\o_{1}$. Since $\o_{3}^{\o}<2^{\o_{3}}$ in the ground
model, in the generic extension, $(2^{\o_{3}})^{V}$ is a cardinal
$>\o_{1}$. So there are $\a <\b <(2^{\o_{3}})^{V}$ such that
$\A_{\a}\not\cong\A_{\b}$ in the extension. Since countable
subsets are not added, E does not have a winning strategy for
$EF^{2}_{\o_{1}}(\A_{\a},\A_{\b})$ (in the ground model). $\eop$

Now Lemmas 2.11 and 2.12 imply Theorem ~\ref{1} (i). $\eop$

Before proving the theorem, we make some remarks which follow from
the proof.

\begin{remark}\label{r}
{\rm In many cases in Theorem~\ref{1}, the assumption on $2^{\o}$
can be removed. For example, this is true of linear orders. An
easy proof for this is given in \cite{Hu}, alternatively this
follows immediately from the proof of Theorem~\ref{1} (i) by
checking where the assumption $2^{\o}<2^{\o_{3}}$ was needed.
Another case where the assumption on $2^{\o}$ can be removed is
the case that
 $\theta =\o_{3}$ in the stable unsuperstable case.
This follows from the proof of Theorem~\ref{1} (iii) by noticing
that the black box can now be replaced by an argument from
\cite{HS}. Another remark is that in Theorem~\ref{1}  (i) and
(ii), $\o_{3}$ can be replaced by any cardinal $\k\ge\o_{3}$ such
that $\k$ is a successor of a regular cardinal and
$2^{\k}>\k^{\o}$. Finally, in Theorem~\ref{1} (iii), $\o_{3}$ can
be replaced by any cardinal $\k\ge\o_{3}$ such that $\k$ is a
successor of a regular cardinal and $\k^{\o}=\k$.}
\end{remark}
}

\end{document}